\documentclass[a4paper,12pt]{article}

\usepackage[dvipdfmx]{graphicx}
\usepackage{amsmath}
\usepackage{url}
\usepackage{txfonts}
\usepackage{multirow}

\begin{document}
\title{A Highly Efficient Implementation of Multiple Precision Sparse Matrix-Vector Multiplication and Its Application to Product-type Krylov Subspace Methods}
\author{Tomonori Kouya\\Shizuoka Institute of Science and Technology\\2200-2 Toyosawa, Fukuroi, Shizuoka 437-8555 Japan}
\maketitle

\begin{quotation}\begin{center}\bf Abstract\end{center}
We evaluate the performance of the Krylov subspace method by using 
highly efficient multiple precision sparse matrix-vector 
multiplication (SpMV). BNCpack is our multiple precision numerical 
computation library based on MPFR/GMP, which is one of the most 
efficient arbitrary precision floating-point arithmetic libraries. 
However, it does not include functions that can manipulate multiple 
precision sparse matrices. Therefore, by using benchmark tests, we 
show that SpMV implemented in these functions can be more efficient. 
Finally, we also show that product-type Krylov subspace methods such 
as BiCG and GPBiCG in which we have embedded SpMV, can efficiently 
solve large-scale linear systems of equations provided in the UF 
sparse matrix collections in a memory-restricted computing 
environment.
\end{quotation}

\section{Introduction}

Large-scale scientific computation is becoming more and more 
necessary these days in a variety of fields. Consequently, a number of 
ill-conditioned problems are created, which cannot be solved precisely 
with standard double precision floating-point arithmetic. To solve 
such ill-conditioned problems, multiple precision floating-point 
arithmetic, such as quadruple or octuple precision is needed as a 
countermeasure.

At present, although high performance numerical software libraries 
that support multiple precision floating-point arithmetic are rare, 
they are being developing. These libraries are built upon primitive 
and well-tuned arithmetic suitable for standard architecture CPUs and 
GPUs. Numerical solutions with the required precision can be obtained 
on a commodity PC.

As a result, increasing more memory is required for multiple precision 
floating-point numbers, and hence, memory needs to be specifically 
allocated for basic linear algebra computation. If the matrix  being 
treated is sparse, the amount of memory can be reduced by storing 
nonzero elements in the matrix. 

Some double precision libraries for treating sparse matrices have been 
provided such as SuiteSparse\cite{suitesparse} and Intel Math Kernel 
(IMK) Library\cite{imkl}, which have direct and iterative methods to 
solve sparse linear systems of equations. If you want to build a 
multiple precision sparse matrix library, you can convert such double 
precision libraries capable of manipulating multiple precision sparse 
matrices, such as MPACK\cite{mpack} which are built to convert LAPACK 
with QD/GD\cite{qd} and GNU MPFR/GMP\cite{mpfr}. However, our multiple 
sparse matrix-vector multiplication is developed as an additional 
function in BNCpack\cite{bnc}, which is our own double and multiple 
precision numerical computation library based on GNU MPFR/GMP, the 
multiple precision floating-point arithmetic library. 

In MPFR/GMP, zero multiplication is made faster to save computational 
time, and we need to estimate how our multiple precision sparse 
matrix-vector product (SpMV for short) can reduce computational time 
in comparison with multiple precision dense matrix-vector 
multiplication (Dense MV) through numerical experiments. The multiple 
precision SpMV is shown to be very efficient and to reduce the 
computational time for product-type Krylov subspace methods in limited 
memory environments.

In this paper, we first introduce the multiple precision SpMV and the 
performance of MPFR/GMP. Next, we benchmark the SpMV in a standard PC 
environment and evaluate the performance in comparison with Dense MV. 
Finally, we apply SpMV to product-type Krylov subspace methods 
provided in BNCpack and demonstrate its efficiency by solving test 
problems using some sparse matrices in the UF Sparse Matrix 
Collection.

\section{Sparse matrix structure based on BNCpack}

In the current BNCpack, all of the multiple precision functions are 
based on MPFR/GMP (MPFR over the MPN kernel of GMP). The significant 
features of MPFR/GMP are as follows:
\begin{itemize}
	\item Arbitrary precision floating-point arithmetic are 
executable according to the IEEE754-1985 standard which supports NaN 
(Not-a-Number), $\pm$Inf (Infinity) and four rounding modes (RN, RZ, 
RP and RM).
	\item The multiple precision floating-point data type 
(\url{mpfr_t}) in MPFR/GMP is defined as a C structure and the 
mantissa part is stored in a separate memory area.
	\item Any mixed precision arithmetic is available everywhere. 
Each of the \url{mpfr_t} variables have their own precision (in bits), 
whereby the precision can change anytime and everywhere throughout the 
running process.
	\item An original acceleration mechanism is embedded to save 
computational time. For example, the arithmetic with fixed precision 
data types such as integer, float and double variables or zero, is 
very fast.
\end{itemize}

\begin{figure}[htbp]
\begin{center}
\includegraphics[width=.8\textwidth]{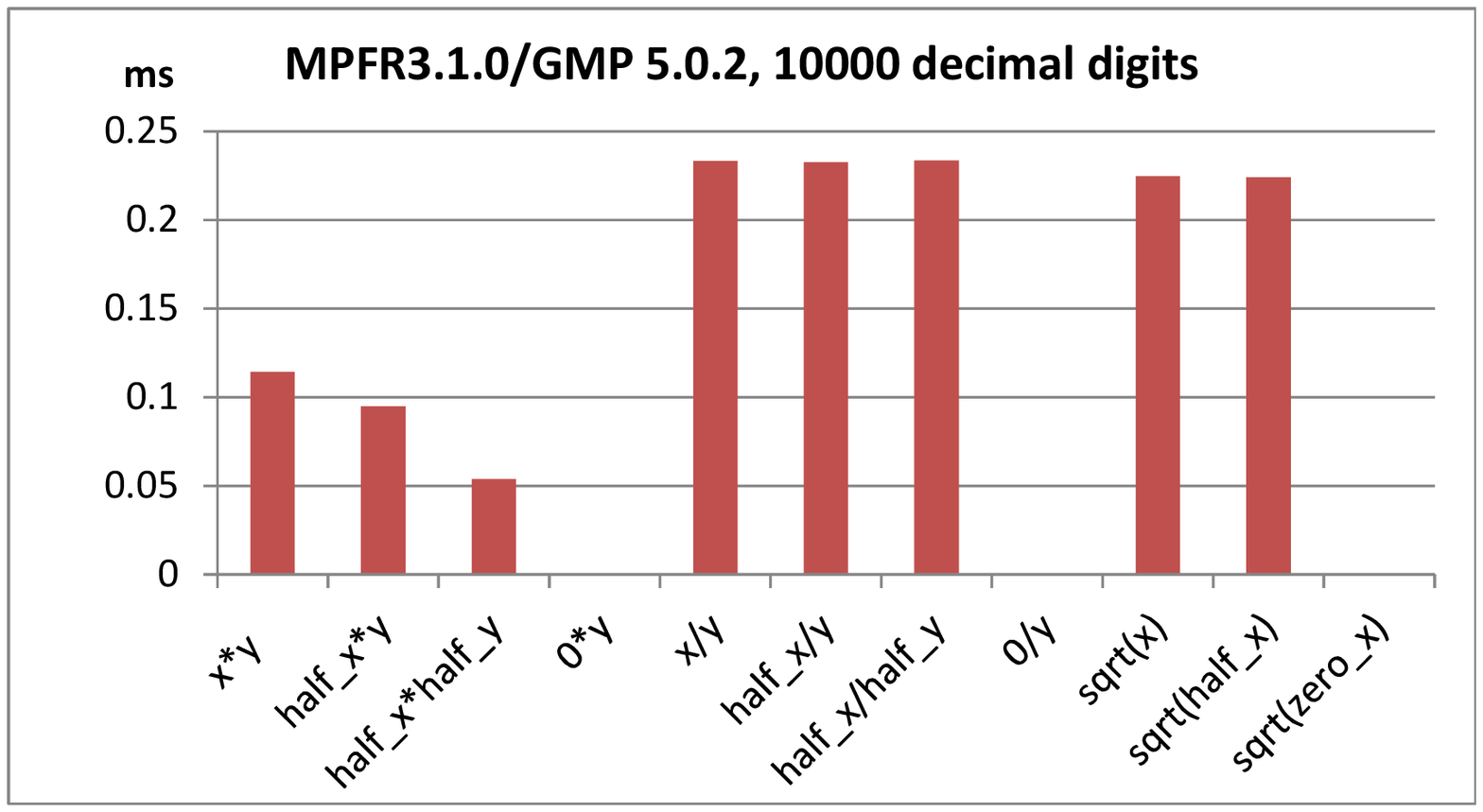}
\caption{Multiplication performance of variables of set precision: MPFR/GMP}\label{fig:mpfr_half} 
\end{center}
\end{figure}

In particular, the original acceleration mechanism of MPFR/GMP is very 
efficient for practical use.  \figurename\ \ref{fig:mpfr_half}) shows 
the performance of multiplication with GMP and MPFR/GMP. The variables 
\verb|x| and \verb|y| have 10000 decimal digits (333220 bits) for the 
mantissa, \verb|half_x| and \verb|half_y| are half digits, and 0 means 
storing zero. Comparing \verb|x*y|, \verb|half_x*half_y| and 
\verb|0*x|, we see that the \verb|0*x| was completed fastest.

Therefore, the particular efficiency of MV can be obtained using a 
dense matrix filled with zero elements. The question then is, what is 
performance of SpMV upon implementing multiple precision sparse matrix 
data types?

Next, we discuss the multiple precision sparse matrix structure based 
on MPFR/GMP. The Compressed Sparse Row (CSR) format was selected as it 
just as the nonzero elements is expressed as a one-dimensional array. 
The C structure is as follows:
\begin{verbatim}
typedef struct {
  unsigned long prec;        // Precision in bits
  mpf_t *element;            // Nonzero elements of the sparse matrix
  long int row_dim, col_dim; // Dimension of row and column
  long int **nzero_index;    // Index of non zero elements
  long int *nzero_col_dim;   // Number of nonzero elements per row
  long int *nzero_row_dim;   // Number of nonzero elements per column
  long int nzero_total_num;  // Total number of nonzero elements
} mpfrsmatrix;
\end{verbatim}

\begin{figure}[htbp]
\begin{center}
\includegraphics[width=.5\textwidth]{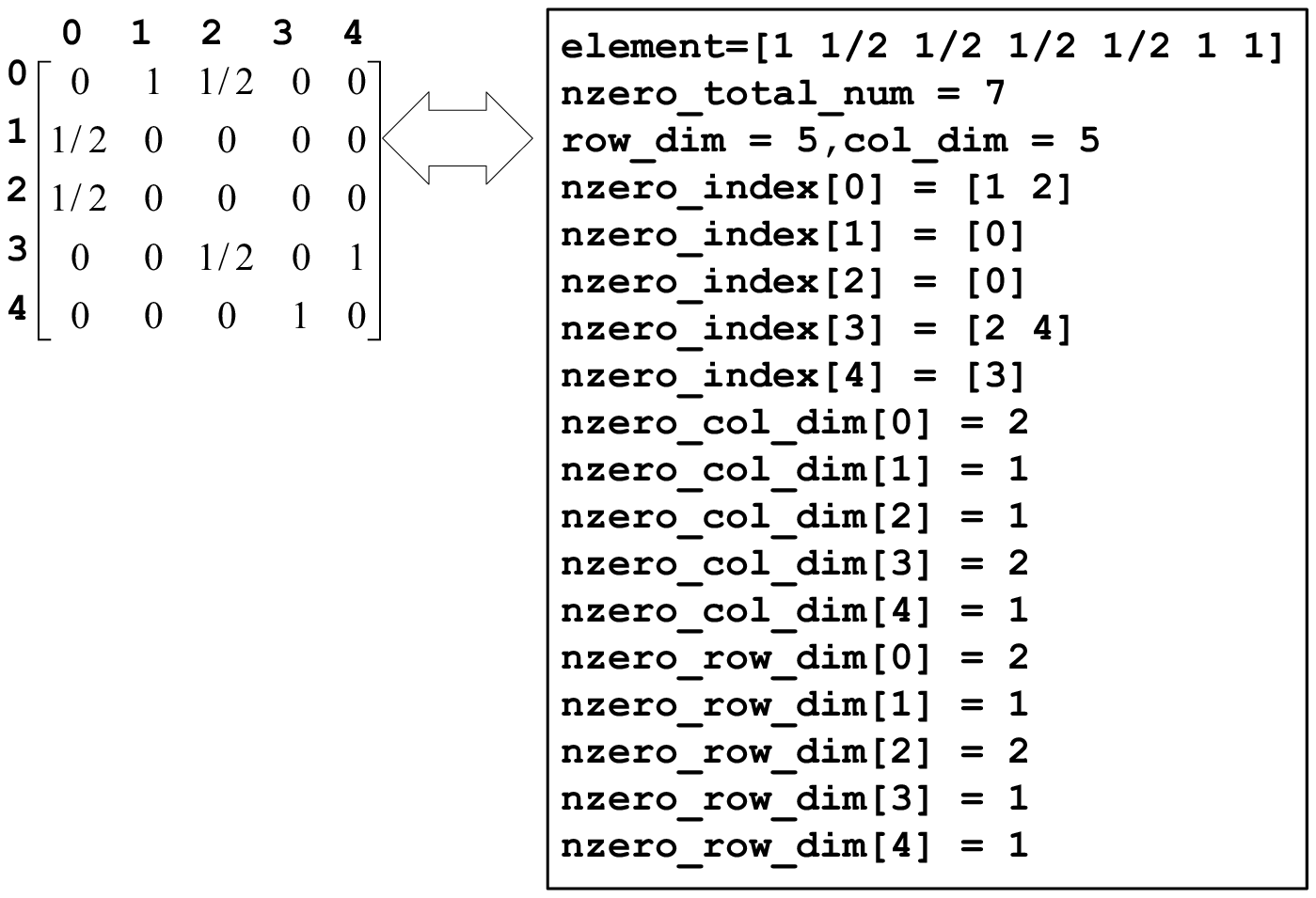}
\caption{Sparse matrix structure implemented in BNCpack.}
\label{fig:spmat_structure}
\end{center}
\end{figure}

According to this structure, for example, the $5\times 5$ random 
sparse matrix can be as shown in \figurename\ 
\ref{fig:spmat_structure}.

As the dimension of matrix increases, the memory size can be reduced 
by using this sparse matrix structure, which can then allow more 
large-scale matrices to be manipulated, which cannot be stored as 
dense matrices. As mentioned above, due to the original acceleration 
mechanism embedded in MPFR/GMP, it is not known to what extent the 
computational time of sparse matrix-vector multiplication is reduced 
without benchmark tests. 

For more details, the source codes used in this paper can be found in 
the latest BNCpack\cite{bnc} library.

%
\section{Performance evaluation of multiple precision SpMV}

We evaluate the performance testing described in this paper with the 
following. Although the environment used has a multi-core CPU, all 
tests are executed in a single thread.
\begin{description}
\item[H/W] Intel Core i7 920, 8GB RAM
\item[S/W] CentOS 5.6 x86\_64, gcc 4.1.2
\item[MPFR/GMP] MPFR 3.0.1, GMP 5.0.1
\item[BNCpack] BNCpack 0.7
\item[BLAS/SparseBLAS] Intel Math Kernel Library 10.3 Update 1
\end{description}

The sparse matrix for evaluation is the Lotkin matrix which is 
randomly filled with zeros.
\[ A = \left[\begin{array}{cccc}
	1 & 1 & \cdots & 1 \\
	1/2 & 1/3 & \cdots & 1/n \\
	\vdots & \vdots & & \vdots \\
	1/n & 1/(n+1) & \cdots & 1/(2n-1)
\end{array}\right] \]
Each element has a mantissa of $p$ precision. The final sparse matrix 
$A^{(s)}$ has $s$\% zero elements. We select a real vector $
\mathbf{b}$ given by
	\[\mathbf{b} = A^{(s, p)}[1\ 2\ \cdots n]^T\]
and evaluate the real matrix-vector multiplication $A^{(s,p)}
\mathbf{b}$. We compare this with $A^{(s)}_{Dense}$, the dense matrix 
format generated by storing the sparse matrix $A^{(s)}$, by examining 
the ratio of the computational time for $A^{(s)}_{Dense}\mathbf{b}$, 
which is defined as the ``Speedup Ratio."

\begin{figure}[htbp]
\begin{center}
\includegraphics[width=.6\textwidth]{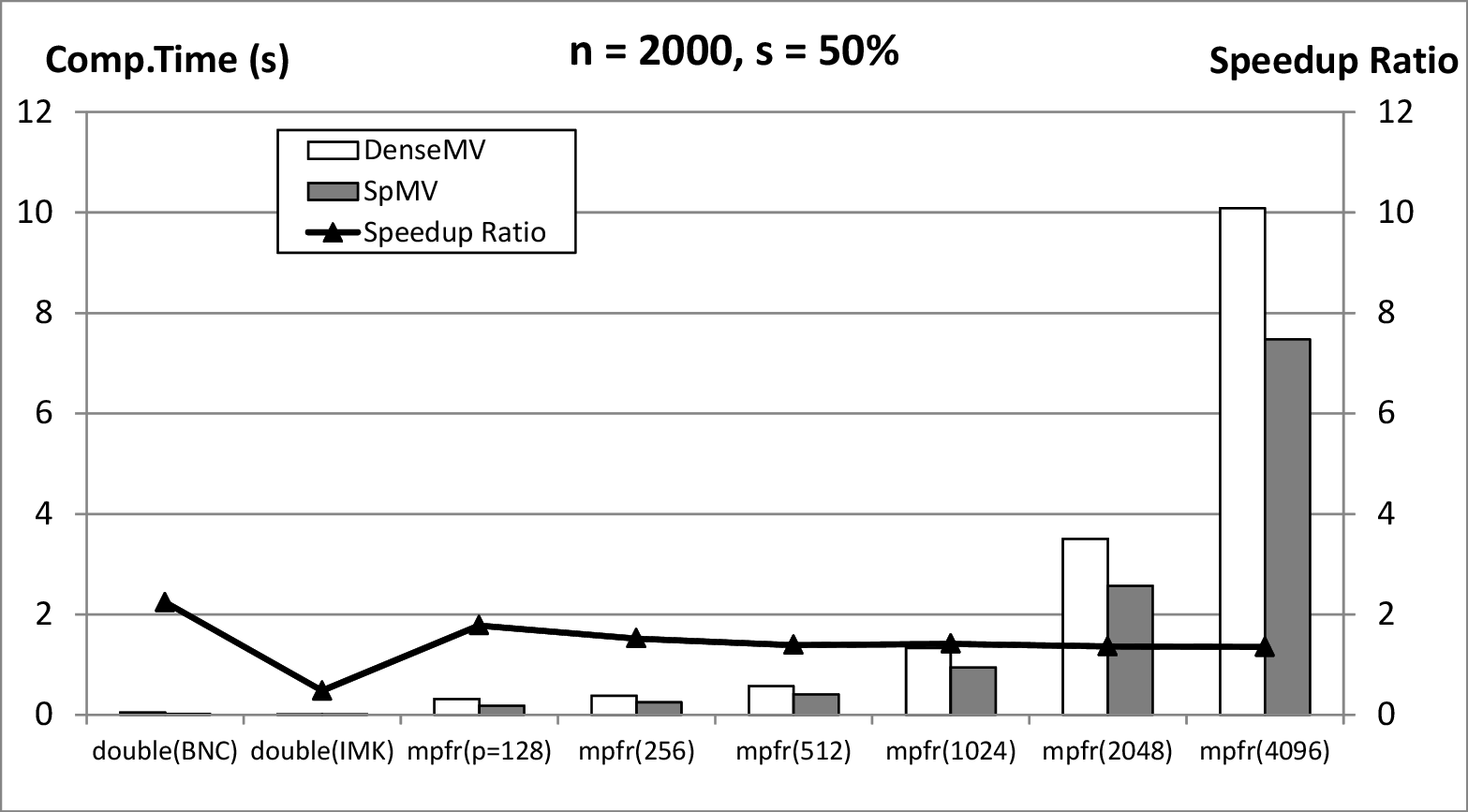}
\includegraphics[width=.6\textwidth]{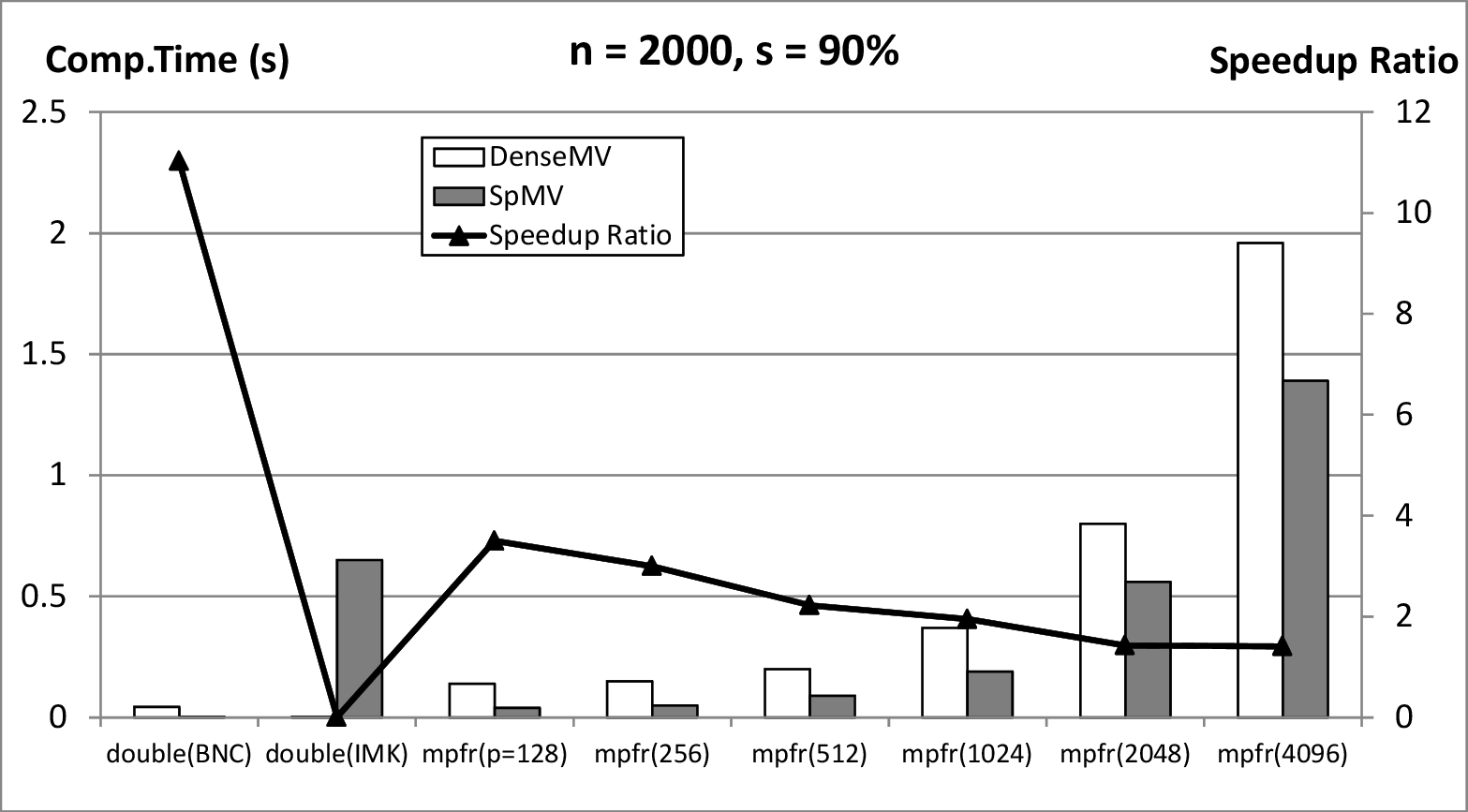}
\caption{Speedup ratio for SpMV (vs. Dense MV).}
\label{fig:sp_benchmark}
\end{center}
\end{figure}

The computational time and the speedup ratio are shown in \figurename\ 
\ref{fig:sp_benchmark}. For reference, the results for the double 
precision BNCpack and Intel Math Kernel (DGEMV/DCOOGEMV) are also 
shown.

\begin{figure}[htbp]
\begin{center}
\includegraphics[width=.95\textwidth]{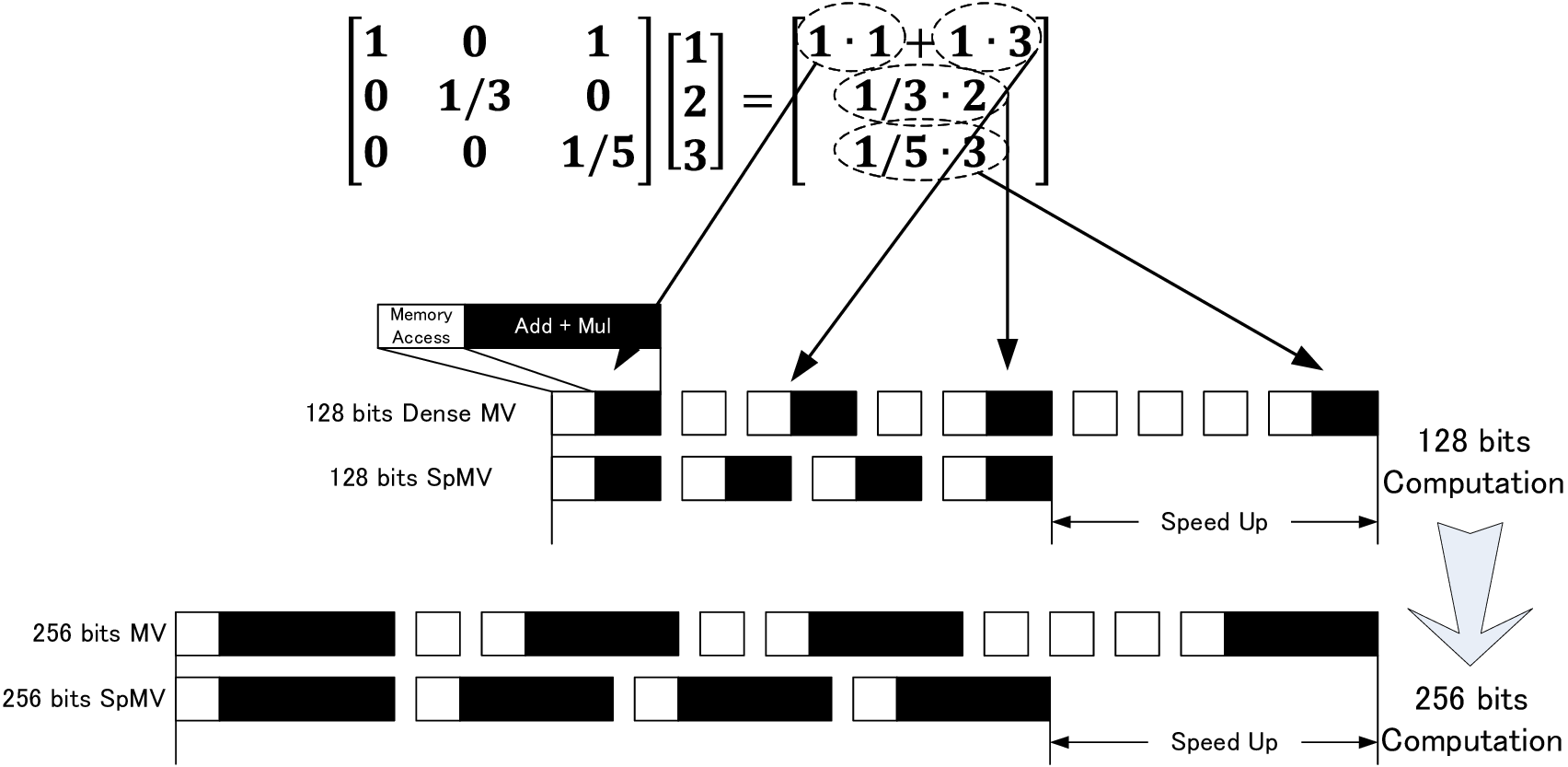}
\caption{Explanation of speedup using SpMV.}
\label{fig:explanation_speedup}
\end{center}
\end{figure}

Double precision SpMV can be effective at 2000-dimesional matrix-
vector multiplication. Multiple precision SpMV can be less effective 
than double precision. The speedup ratio varies between about 1.3 and 
3.2. The reason can be interpreted to be the effect of the MPFR/GMP 
acceleration mechanism, whereby zero multiplication can be executed 
very quickly in dense MV. This can be understood from the fact that 
the higher the precision in SpMV, the lower is the speedup ratio 
actually gained. \figurename\ \ref{fig:explanation_speedup} explains 
this mechanism.

The performance evaluation of SpMV shows that the computational time 
can be reduced if $s$ and $p$ are sufficiently large and that the 
speedup of SpMV is effective if the matrix is large and the precision 
is not very high.

\section{Application to product-type Krylov subspace methods}

The performance evaluation of SpMV described above demonstrates that 
not only is memory saved but speedup is also achieved. Consequently 
more efficient product-type Krylov subspace methods\cite{Zhang2002297} 
with multiple precision SpMV can be implemented. We selected the BiCG, 
CGS, BiCGSTAB and GPBiCG methods.

Krylov subspace methods and their variations are well-known for their 
sensitivity to the effect of round-off errors occurring during their 
iterative processes, and this effect can be reduced by using multiple 
precision floating-point arithmetic\cite{Saito2011}. In this section, 
we describe the performance evaluation and numerical experiments using 
multiple precision product-type Krylov subspace methods. We also 
describe the numerical experiments with the ILU(0) method for left 
preconditioning.

%
\subsection{Performance evaluation of multiple precision product-type 
Krylov subspace methods}

From the implemented BiCG, CGS, BiCGSTAB, and GPBiCG, we present the 
algorithms of BiCG and GPBiCG (\figurename\ \ref{fig:gpbicg}. The 
underlined sections indicate the replacement for SpMV.

%
\begin{figure}[htbp]
\begin{center}
\begin{minipage}{.42\textwidth}\small
\begin{description}
\item[] $\mathbf{x}_0$: Initial guess
\item[] $\mathbf{r}_0$: Initial residual ($\mathbf{r}_0 = \mathbf{b} - 
A\mathbf{x}_0$)
\item[] $\widetilde{\mathbf{r}_0}$: $= \mathbf{r}_0$.
\item[] $K$: Matrix for Left preconditioning.
\item[] {\bf for} $i = 1, 2, ...$
\begin{description}
	\item[] Solve $K\mathbf{w}_{i-1} = \mathbf{r}_{i-1}$ on $
\mathbf{w}_{i-1}$.
	\item[] Solve $K^T\ \widetilde{\mathbf{w}_{i-1}}= 
\widetilde{\mathbf{r}_{i-1}}$ on $\widetilde{\mathbf{w}_{i-1}}$.
	\item[] $\rho_{i-1} = (\widetilde{\mathbf{w}_{i-1}}, 
\mathbf{w}_{i-1})$
	\item[] {\bf if} $\rho_{i-1}=0$ {\bf then} exits.
	\item[] {\bf if} $i = 1$ {\bf then}\strut\par
	\begin{itemize}
		\item[] $\mathbf{p}_1 = \mathbf{w}_0$
		\item[] $\widetilde{\mathbf{p}_1} = 
\widetilde{\mathbf{w}_0}$
	\end{itemize}
	\item[] {\bf else}\strut\par
	\begin{itemize}
		\item[] $\beta_{i-1} = \rho_{i-1} / \rho_{i-2}$
		\item[] $\mathbf{p}_i = \mathbf{w}_{i-1} + \beta_{i-
1}\mathbf{p}_{i-1}$
		\item[] $\widetilde{\mathbf{p}_i} = 
\widetilde{\mathbf{w}_i} + \beta_{i-1}\widetilde{\mathbf{p}_{i-1}}$
	\end{itemize}
	\item[] {\bf end if}
	\item[] \underline{$\mathbf{z}_i = A\mathbf{p}_i$}
	\item[] \underline{$\widetilde{\mathbf{z}_i} = A
\widetilde{\mathbf{p}_i}$}
	\item[] $\alpha_i = \rho_{i-1} / (\widetilde{\mathbf{p}_i}, 
\mathbf{z}_i)$
	\item[] $\mathbf{x}_i = \mathbf{x}_{i-1} + \alpha_i 
\mathbf{p}_i$
	\item[] $\mathbf{r}_i = \mathbf{r}_{i-1} - \alpha_i 
\mathbf{z}_i$
	\item[] $\widetilde{\mathbf{r}_i} = \widetilde{\mathbf{r}_{i-
1}} - \alpha_i \widetilde{\mathbf{z}_i}$
	\item[] Convergence check(*).
\end{description}
\item[] {\bf end for}
\end{description}
\end{minipage}
%
%
\begin{minipage}{.55\textwidth}\small
\begin{description}
\item[] $\mathbf{u} = \mathbf{z} = 0$
\item[] {\bf for} $i = 1, 2, ...$
\begin{description}
	\item[] $\rho_{i-1} = (\widetilde{\mathbf{r}}, \mathbf{r}_{i-
1})$
	\item[] {\bf if} $\rho_{i-1}=0$ {\bf then} exits.
	\item[] {\bf if} $i = 1$ {\bf then}\strut\par
	\begin{itemize}
		\item[] $\mathbf{p} = \mathbf{r}_0$
		\item[] \underline{$\mathbf{q} = A\mathbf{p}$}
		\item[] $\alpha_i = \rho_{i-1} / (\widetilde{r}, 
\mathbf{q})$
		\item[] $\mathbf{t} = \mathbf{r}_{i-1} - \alpha_i 
\mathbf{q}$
		\item[] \underline{$\mathbf{v} = A\mathbf{t}$}
		\item[] $\mathbf{y} = \alpha_i \mathbf{q} - 
\mathbf{r}_{i-1}$
		\item[] $\mu_2 = (\mathbf{v}, \mathbf{t})$, $\mu_5 = 
(\mathbf{v}, \mathbf{v})$
		\item[] $\zeta = \mu_2 / \mu_5$
		\item[] $\eta = 0$
	\end{itemize}
	\item[] {\bf else}\strut\par
	\begin{itemize}
		\item[] $\beta_{i-1} = (\rho_{i-1} / \rho_{i-2})
(\alpha_{i-1}/\zeta)$
		\item[] $\mathbf{w} = \mathbf{v} + \beta_{i-
1}\mathbf{q}$
		\item[] $\mathbf{p} = \mathbf{r}_{i-1} + \beta_{i-1}
(\mathbf{p} - \mathbf{u})$
		\item[] \underline{$\mathbf{q} = A\mathbf{p}$}
		\item[] $\alpha_i = \rho_{i-1} / (\widetilde{r}, 
\mathbf{q})$
		\item[] $\mathbf{s} = \mathbf{t} - \mathbf{r}_{i-1}$
		\item[] $\mathbf{t} = \mathbf{r}_{i-1} - \alpha_i 
\mathbf{q}$
		\item[] \underline{$\mathbf{v} = A\mathbf{t}$}
		\item[] $\mathbf{y} = \mathbf{s} - \alpha_i 
(\mathbf{w} - \mathbf{q})$
		\item[] $\mu_1 = (\mathbf{y}, \mathbf{y})$; $\mu_2 = 
(\mathbf{v}, \mathbf{t})$
		\item[] $\mu_3 = (\mathbf{y}, \mathbf{t})$; $\mu_4 = 
(\mathbf{v}, \mathbf{y})$
		\item[] $\mu_5 = (\mathbf{v}, \mathbf{v})$;
		\item[] $\tau = \mu_5 \mu_1 - \overline{\mu_4}\mu_4$
		\item[] $\zeta = (\mu_1 \mu_2 - \mu_3 \mu_4) / \tau$
		\item[] $\eta = (\mu_5 \mu_3 - \overline{\mu_4} \mu_2) 
/ \tau$
	\end{itemize}
	\item[] {\bf end if}
	\item[] $\mathbf{u} = \zeta\mathbf{q} + \eta(\mathbf{s} + 
\beta_{i-1}\mathbf{u})$
	\item[] $\mathbf{z} = \zeta\mathbf{r}_{i-1} + \zeta\mathbf{z} 
- \alpha_i\mathbf{u}$
	\item[] $\mathbf{x}_i = \mathbf{x}_{i-1} + \alpha_i \mathbf{p} 
+ \mathbf{z}$
	\item[] Convergence check(*)
	\item[] $\mathbf{r}_i = \mathbf{t} - \eta \mathbf{y} - \zeta
\mathbf{u}$
	\item[] exits if $\zeta = 0$.
\end{description}
\item[] {\bf end for}
\item[]
\end{description}
\end{minipage}
\end{center}
\caption{Algorithms for BiCG (left) and GPBiCG (right). The underlined 
sections indicate the replacement for SpMV.}\label{fig:gpbicg}
\end{figure}

In the numerical experiments below, all the sparse matrices used are 
provided in the UF Sparse Matrix Collection\cite{ufsparse}. The 
original MTX files storing the sparse matrices have only double 
precision elements and hence, these are converted into multiple 
precision floating-point numbers. The true solutions are commonly $
\mathbf{x} = [1\ 2\ \cdots\ n]^T$, and thus, the constant vectors are 
$\mathbf{b} = A\mathbf{x}$, which are calculated in multiple precision 
computation. The convergence condition ((*) in \figurename\ 
\ref{fig:gpbicg}) is satisfied if \[ ||\mathbf{r}_k ||_2 \leq 10^{-20} 
|| \mathbf{r}_0 ||_2 + 10^{-50}. \]

\paragraph{DRICAV/cavity04}

DRIVCAV/cavity04 ($n = 317$) is selected for small-scale problems. The 
number of nonzero elements is 7327 ($s = 93.7$\%). Problems of the 
size of this scale can be completely stored in dense matrix format in 
RAM, which enables us to find the speedup ratio for the usage of SpMV. 
The precision of the elements of the matrix is extended from 512 bits 
(about 154 decimal digits) to 8192 bits (about 2466 decimal digits) 
and we compare the speedup ratio with Dense MVs.

\begin{figure}[htb]
\begin{center}
\includegraphics[width=.95\textwidth]{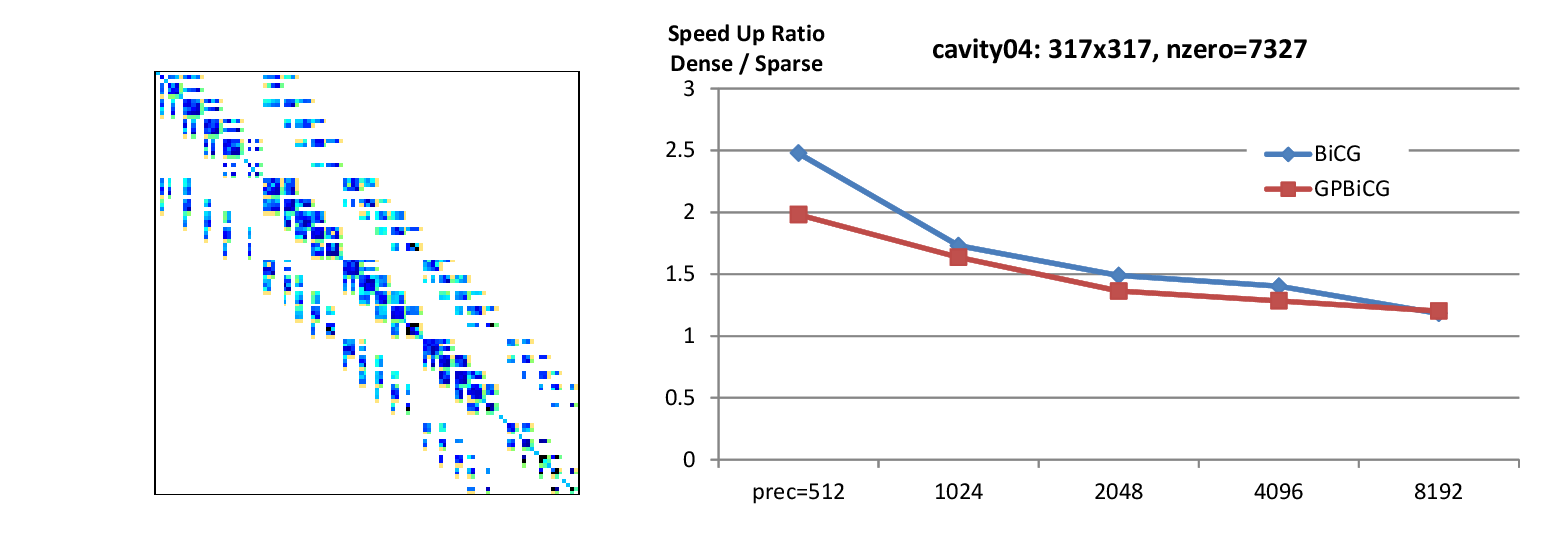}
\caption{Matrix structure of cavity04 (left) and speedup ratio 
(right).}\label{fig:cavity04}
\end{center}
\end{figure}

The structure of cavity04 matrix and the speedup ratio for BiCG and 
GPBiCG are shown in \figurename\ \ref{fig:cavity04}.

The iterative and computational times for the BiCG, CGS, BiCGSTAB, and 
GPBiCG methods for the cavity04 problem are listed in \tablename\ 
\ref{table:cavity04_iterative_time}.

\begin{table}
\begin{center}
\caption{Iterative and Computational Times (s): 
cavity04}\label{table:cavity04_iterative_time}
\begin{tabular}{|c|c|c|c|c|c|}\hline
Iterative Times (s)	& prec=512	& 1024	& 2048	& 4096	& 8192 
\\ \hline
BiCG			& Not Convergent & NC & 236 (4.83) & 236 (13.74) & 236 (39.91) \\ \hline
CGS			& NC & NC & NC & 236 (13.74)	& 236 (40.29) \\ \hline
BiCGSTAB		& NC & NC & NC & 260 (16.23)	& 236 (42.79) \\ \hline
GPBiCG		& NC & NC & NC & 24 3(23.99)	& 236 (67.48) \\ 
\hline
\end{tabular}
\end{center}
\end{table}

\begin{figure}[htb]
\begin{center}
\includegraphics[width=.7\textwidth]{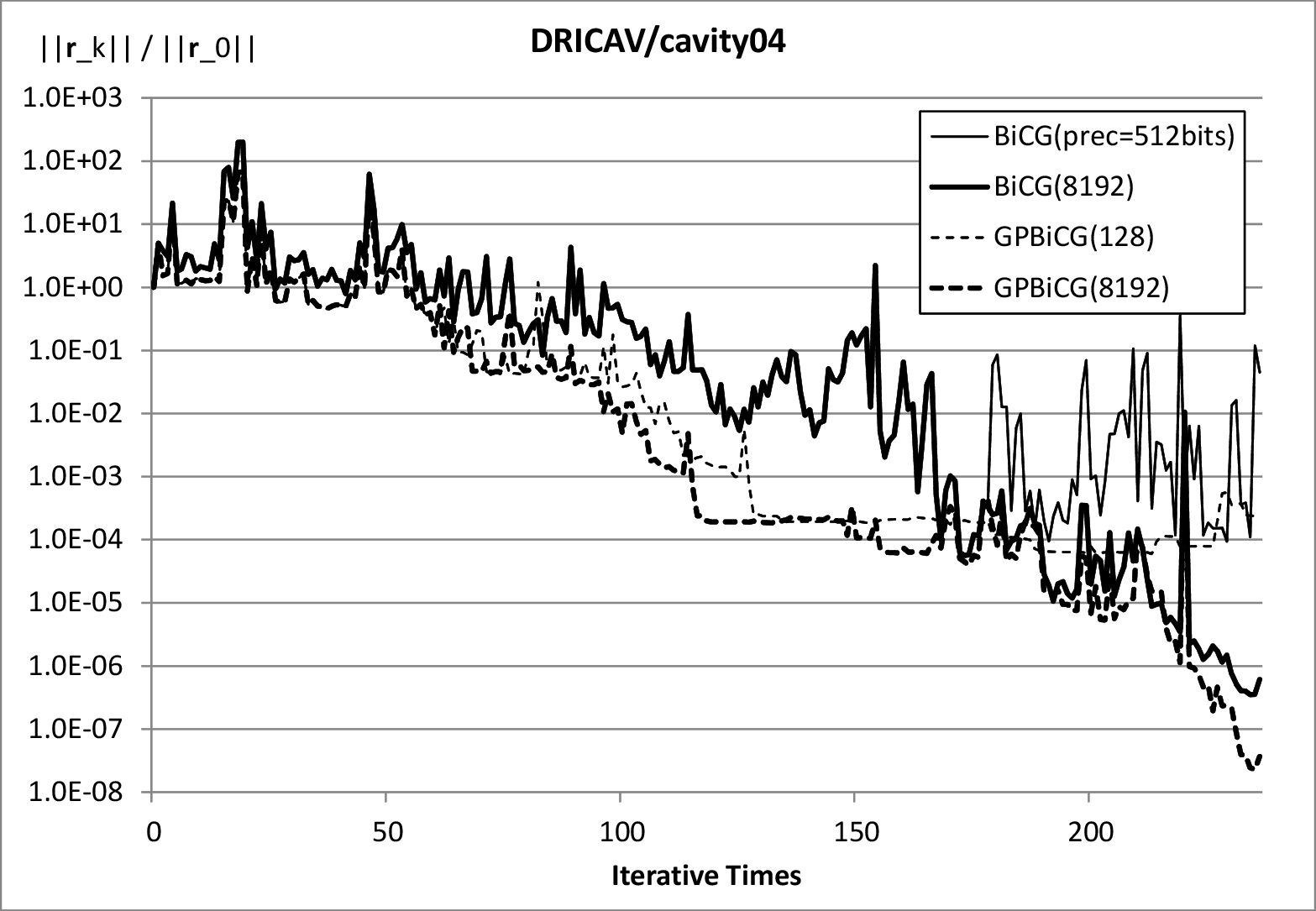}
\caption{History of $||\mathbf{r}_k||_2 / ||\mathbf{r}_0||_2$ for the 
cavity04 problem}\label{fig:cavity04_res_history}
\end{center}
\end{figure}

In the case of the cavity04 problem, the history of residual norms 
(\figurename\ \ref{fig:cavity04_res_history}) shows that product-type 
Krylov subspace methods need more than 2048 bits precision to 
converge. Multiple precision arithmetic is very effective in obtaining 
more precise numerical solutions.

\paragraph{Averous/epb3}

Averous/epb3 ($n=84617$) is selected for larger scale problems. The 
true solution $\mathbf{x}$ and the constant vector $\mathbf{b}$ are 
the same as for the cavity04 problem.

The number of nonzero elements is 463,625 ($s = 99.0035$\%), and 
hence, if the matrix is stored in dense matrix format, the total 
amount of memory required is about 320 GB in the case of 128-bit 
precision and about 70 TB in the case of 8192-bit precision. With our 
sparse matrix format, about 3.5 MB is required for 128-bit precision 
and about 456 MB for 8192-bit precision. Hence, numerical experiments 
with the epb3 matrix can be carried out on standard memory-restricted 
PCs.

The iterative and computational times are listed in \tablename\ 
\ref{table:epb3}. ``NC" means that the process cannot converge 
properly.

\begin{table}
\begin{center}
\caption{Iterative and computational (in parenthesis) times for the 
epb3 problem.}\label{table:epb3}
\begin{tabular}{|c|c|c|c|c|c|} \hline
       & 512bits & 1024 & 2048           & 4096  & 8192 \\ \hline
BiCG   & NC      & NC   & 8351(4335 s) & 5251(7824 s) & 3567(15551 s) \\ \hline
GPBiCG & NC      & NC   & NC             & 7039(12473 s) & 4449(24179 s) \\ \hline
\end{tabular}
\end{center}
\end{table}

In the epb3 problem, the effectiveness of multiple precision arithmetic 
is proven because higher precision can reduce the iterative times of 
every product-type Krylov subspace methods.

%
\subsection{Performance evaluation of the BiCGSTAB method with ILU(0) 
for left preconditioning}

For the double precision Krylov subspace method and its variations, we 
normally set up the matrix $K$ as a precondition for each problems in 
order to speed up the iteration. In this case, we use the BiCGSTAB 
method with incomplete LU decomposition (ILU(0)) as the left 
precondition, and examine the numerical properties and efficiency of 
ILU(0).

For the cavity04 problem, the results of BiCGSTAB with ILU(0) are 
summarized in \tablename\ \ref{table:bicgstab_ilu}.

\begin{table}[htbp]
\begin{center}
\caption{BiCGSTAB method with ILU(0) for left preconditioning: 
cavity04}\label{table:bicgstab_ilu}
\begin{tabular}{|c|c|c|c|c|c|c|}\hline
\phantom{cavity04}	& \phantom{dim=317}	& prec=512	& 1024	& 2048	& 4096	& 8192 \\ \hline
\multirow{2}{9em}{BiCGSTAB}	& Iter. Times	& NC	& NC	& NC	& 260	& 236 \\ \cline{2-7}
		& Time (s)	& \phantom{4.25}& \phantom{5.55}	& \phantom{8.26}	& \underline{16.23}	& 42.79 \\ \hline
\multirow{2}{9em}{BiCGSTAB +ILU(0)} & Iter. Times	& NC	& NC	& 274	& 238	& 234 \\ \cline{2-7}
				 & Time (s)	& \phantom{81.65}	& \phantom{53.8}	& \underline{45.94}	& 58.75	& 109.29 \\ \hline
\end{tabular}
\end{center}
\end{table}

In this case, 4096-bit precision computation is needed for convergence 
for the simple BiCGSTAB method. Although the application of ILU(0) can 
result in convergence of the BiCGSTAB process, the simple 4096-bit 
BiCGSTAB method is more efficient because the computational time with 
ILU(0) is long. The results  of the numerical experiments show that 
preconditioning in multiple precision computation is not efficient due 
to the effect of the matrix structure and other such factors, if it 
performs better in double precision computation.

\section{Conclusion}

Although multiple precision computation is needed in various types of 
scientific computation, software libraries applicable to large-scale 
problems are very rare. Our multiple precision SpMV, based on BNCpack 
and MPFR/GMP and the product-type Krylov subspace methods, is 
particularly useful. The numerical in this paper reveals that multiple 
precision SpMV can save memory space in RAM and can also reduce 
computational time.

In a future work, we will apply these implemented Krylov subspace 
methods to large-scale scientific computation in particular, that 
which includes ordinary and partial differential equations. We will 
also attempt to further validate BNCpack with SpMV and Krylov subspace 
methods by applying it to actual problems.


\end{document}